# THE DIMENSION OF THE SLE CURVES

### By Vincent Beffara

*CNRS—UMPA—ENS Lyon*

Let $\gamma$ be the curve generating a Schramm–Loewner Evolution (SLE) process, with parameter $\kappa \geq 0$. We prove that, with probability one, the Hausdorff dimension of $\gamma$ is equal to $\mathrm{Min}(2, 1 + \kappa/8)$.

**Introduction.** It has been conjectured by theoretical physicists that various lattice models in statistical physics (such as percolation, Potts model, Ising model, uniform spanning trees), taken at their critical point, have a continuous conformally invariant scaling limit when the mesh of the lattice tends to 0. Recently, Oded Schramm [15] introduced a family of random processes which he called Stochastic Loewner Evolutions (or SLE), that are the only possible conformally invariant scaling limits of random cluster interfaces (which are very closely related to all above-mentioned models).

An SLE process is defined using the usual Loewner equation, where the driving function is a time-changed Brownian motion. More specifically, in the present paper we will be mainly concerned with SLE in the upper-half plane (sometimes called chordal SLE), defined by the following PDE:

$$(0.1) \qquad \partial_t g_t(z) = \frac{2}{g_t(z) - \sqrt{\kappa} B_t}, \qquad g_0(z) = z,$$

where $(B_t)$ is a standard Brownian motion on the real line and $\kappa$ is a positive parameter. It can be shown that this equation defines a family $(g_t)$ of conformal mappings from simply connected domains $(H_t)$ contained in the upper-half plane, onto $\mathbb{H}$. We shall denote by $K_t$ the closure of the complement of $H_t$ in $\mathbb{H}$: then for all $t > 0$, $K_t$ is a compact subset of $\overline{\mathbb{H}}$ and the family $(K_t)$ is increasing. For each value $\kappa > 0$, this defines a random process denoted by $\mathrm{SLE}_\kappa$ (see, e.g., [14] for more details on SLE).

There are very few cases where convergence of a discrete model to $\mathrm{SLE}_\kappa$ is proved: Smirnov [17] (see also the related work of Camia and Newman [3]) showed that SLE₆ is the scaling limit of critical site percolation interfaces on the triangular grid, and Lawler, Schramm and Werner [12] have proved that $\text{SLE}_2$ and $\text{SLE}_8$ are the respective scaling limits of planar loop-erased random walks and uniform Peano curves. Convergence of the "harmonic explorer" was obtained by Schramm and Sheffield [16], and there is also strong evidence [13] that the infinite self-avoiding walk in the half-plane is related to $\text{SLE}_{8/3}$.

It is natural to study the geometry of $\text{SLE}_\kappa$, and in particular, its dependence on $\kappa$. It is known (cf. [12, 14]) that, for each $\kappa > 0$, the process $(K_t)$ is generated by a random curve $\gamma : [0, \infty) \to \overline{\mathbb{H}}$ (called the *trace* of the SLE or the SLE curve), in the following sense: For each $t > 0$, $H_t$ is the unique unbounded connected component of $\mathbb{H} \setminus \gamma([0, t])$. Furthermore (see [14]), $\gamma$ is a simple curve when $\kappa \leq 4$, and it is a space-filling curve when $\kappa \geq 8$. The geometry of this curve will be our main object of interest in the present paper.

It is possible, for each $z \in \mathbb{H}$, to evaluate the asymptotics when $\varepsilon \to 0$ of the probability that $\gamma$ intersects the disk of radius $\varepsilon$ around $z$. When $\kappa < 8$, this probability decays like $\varepsilon^s$ for some $s = s(\kappa) > 0$. This (loosely speaking) shows that the expected number of balls of radius $\varepsilon$ needed to cover $\gamma([0, 1])$ (say) is of the order of $\varepsilon^{s-2}$, and implies that the Hausdorff dimension of $\gamma$ is not larger than $2 - s$. Rohde and Schramm [14] used this strategy to show that almost surely the Hausdorff dimension of the $\text{SLE}_\kappa$ trace is not larger than $1 + \kappa/8$ when $\kappa \leq 8$, and they conjectured that this bound was sharp.

Our main result in the present paper is the proof of this conjecture, namely:

THEOREM.  *Let* $(K_t)$ *be an* $\text{SLE}_\kappa$ *in the upper-half plane with* $\kappa > 0$, *let* $\gamma$ *be its trace and let* $\mathcal{H} := \gamma([0, \infty))$. *Then, almost surely,*

$$\dim_H(\mathcal{H}) = 2 \wedge \left(1 + \frac{\kappa}{8}\right).$$

This result was known for $\kappa \geq 8$ (because the curve is then space-filling), $\kappa = 6$ (see [2], recall that this corresponds to the scaling limit of critical percolation clusters) and $\kappa = 8/3$ (this follows from the description of the outer frontier of $\text{SLE}_6$—or planar Brownian motion—in terms of $\text{SLE}_{8/3}$ in [11], and the determination of the dimension of this boundary, see [8, 9]). Note that in both these special cases, the models have a lot of independence built in (the Markov property of planar Brownian motion, the locality property of $\text{SLE}_6$), and that the proofs use it in a fundamental way.

$\text{SLE}_2$ is the scaling limit (see [12]) of the two-dimensional loop-erased walk: Hence, we prove that the Hausdorff dimension of this scaling limit is $5/4$, that is, it is equal to the *growth exponent* of the loop-erased walk (obtained by Kenyon, cf. [6]) and, at least heuristically, this is not surprising. It is not known whether Kenyon's result can be derived using SLE methods.

This exponent $s$ and various other exponents describing exceptional subsets of $\gamma$ are closely related to critical exponents that describe the behavior near the critical point of some functionals of the corresponding statistical physics model. The value of the exponents $1 + \kappa/8$ appear in the theoretical physics literature (see, e.g., [4] for a derivation based on quantum gravity, and the references therein) in terms of the central charge of the model. Let us stress that in the physics literature, the derivation of the exponent is often announced in terms of (almost sure) fractal dimension, thereby omitting to prove the lower bound on the dimension. It might a priori be the case that the value $\varepsilon^{s-2}$ is due to exceptional realizations of $SLE_\kappa$ with exceptionally many visited balls of radius $\varepsilon$, while "typical" realizations of $SLE_\kappa$ meet far fewer disks, in which case the dimension of the curve could be smaller than $2 - s$.

One standard way to exclude such a possibility and to prove that $2 - s$ corresponds to the almost sure dimension of a random fractal is to estimate the variance of the number of $\varepsilon$-disks needed to cover it. This amounts to computing second moments, that is, given *two* balls of radius $\varepsilon$, to estimating the probability that the SLE trace intersects both of them—and this is the hairy part of the proof, especially if there is a long-range dependence in the model. One also needs another nontrivial ingredient: One has to evaluate precisely (i.e., up to multiplicative constants) the probability of intersecting one ball. Even in the Brownian case (see, e.g., [10]), this is not an easy task.

Note that the discrete counterpart of our theorem in the cases $\kappa = 6$ and $\kappa = 2$ is still an open problem. It is known that for critical percolation interfaces (see [18]) and for loop-erased random walks [6], the expected number of steps grows in the appropriate way when the mesh of the lattice goes to zero, but its almost sure behavior is not yet well-understood: For critical percolation, the up-to-constant estimate of the first moment is missing, and for loop-erased random walks, we lack the second moment estimate.

Another natural object is the *boundary* of an SLE, namely, $\partial K_t \cap \mathbb{H}$. For $\kappa \leq 4$, since $\gamma$ is a simple curve, the boundary of the SLE is the SLE itself; for $\kappa > 4$, it is a strict subset of the trace, and it is conjectured to be closely related to the curve of an $SLE_{16/\kappa}$ (this is called *SLE duality*)—in particular, it should have dimension $1 + 2/\kappa$. Again, the first moment estimate is known for all $\kappa$ (though not up to constants), and yields the upper bound on the dimension. The lower bound is known to hold for $\kappa = 6$ (see [8]). A consequence of our main theorem is that it also holds for $\kappa = 8$, because of the continuous counterpart of the duality between uniform spanning trees and loop-erased random walks (which is the basis of Wilson's algorithm, cf. [19]).

The derivation of the lower bound on the dimension relies on the construction of a random Frostman measure $\mu$ supported on the curve. It appears that the properties of this measure are closely related to some of those exhibited by conformal fields—more specifically, the correlations between the measures of disjoint subsets of $\mathbb{H}$ behave similarly to the (conjectured) correlation functions in conformal field theory. See, for instance, Friedrich and Werner [5].

The plan of this paper is as follows. In the first section we review some facts that can be found in our previous paper ([2]) and that we will need later. Section 2 is devoted to the derivation of the up-to-constants estimate of the first moment of the number of disks needed to cover the curve. In Section 3 we will derive the upper bound on the second moment, which will conclude the proof of the main theorem. In the final sections we will comment on the properties of the Frostman measure supported on the SLE curve and on the dimension of the outer boundary of $SLE_8$.

**1. Preliminaries.** As customary, the Hausdorff dimension of the random fractal curve $\gamma$ will be determined using first and second moments estimates. This framework was also used in [2]. We now briefly recall without proofs some tools from that paper that we will use. The following proposition is the continuous version of a similar discrete construction due to Lawler (cf. [7]).

Let $\lambda$ be the Lebesgue measure in $[0,1]^2$, and $(C_\varepsilon)_{\varepsilon>0}$ be a family of random Borel subsets of the square $[0,1]^2$. Assume that for $\varepsilon < \varepsilon'$ we have almost surely $C_\varepsilon \subseteq C_{\varepsilon'}$, and let $C = \cap C_\varepsilon$. Finally, let $s$ be a nonnegative real number. Introduce the following conditions:

1. For all $x \in [0,1]^2$,
$$P(x \in C_\varepsilon) \asymp \varepsilon^s$$
   (where the symbol $\asymp$ means that the ratio between both sides of the expression is bounded above and below by finite positive constants);

2. There exists $c > 0$ such that, for all $x \in [0,1]^2$ and $\varepsilon > 0$,
$$P(\lambda(C_\varepsilon \cap \mathcal{B}(x,\varepsilon)) > c\varepsilon^2 | x \in C_\varepsilon) \geq c > 0;$$

3. There exists $c > 0$ such that, for all $x, y \in [0,1]^2$ and $\varepsilon > 0$,
$$P(\{x,y\} \subset C_\varepsilon) \leq c\varepsilon^{2s}|x-y|^{-s}.$$

PROPOSITION 1.  *With the previous notation:*

1. *If conditions 1 and 2 hold, then a.s.* $\dim_H(C) \leq d - s$;
2. *If conditions 1 and 3 hold, then with positive probability* $\dim_H(C) \geq d - s$.

REMARK. The similar proposition which can be found in [7] is stated in a discrete setup in which condition 2 does not appear. Indeed, in most cases, this condition is a direct consequence of condition 1 and the definition of $C_\varepsilon$ (e.g., if $C_\varepsilon$ is a union of balls of radius $\varepsilon$ as will be the case here).

The value of the exponent in condition 1 is usually given in terms of the principal eigenvalue of a diffusion generator (cf. [1] for further reference). The rule of thumb is as follows:

LEMMA 2. *Let $(X_t)$ be the diffusion on the interval $[0,1]$ generated by the following stochastic differential equation:*

$$dX_t = \sigma \, dB_t + f(X_t) \, dt,$$

*where $(B_t)$ is a standard real-valued Brownian motion, $\sigma$ is a positive constant, and where $f$ is a smooth function on the open unit interval satisfying suitable conditions near the boundary. Let $L$ be the generator of the diffusion, defined by*

$$L\phi = \frac{\sigma^2}{2}\phi'' + f\phi',$$

*and let $\lambda$ be its leading eigenvalue. Then, as $t$ goes to infinity, the probability $p_t$ that the diffusion is defined up to time $t$ tends to 0 as*

$$p_t \asymp e^{-\lambda t}.$$

We voluntarily do not state the conditions satisfied by $f$ in detail here (roughly, $f$ needs to make both 0 and 1 absorbing boundaries, while being steep enough to allow a spectral gap construction—cf. [2] for a more complete statement), because we shall not use the lemma in this form in the present paper; we include it mainly for background reference.

The next two sections contain derivations of conditions 1 and 3; together with Proposition 1, this implies that

$$P\left[\dim_H \mathcal{H} = 1 + \frac{\kappa}{8}\right] > 0.$$

The main theorem then follows from the zero-one law derived in [2], namely:

LEMMA 3 (0–1 law for the trace). *For all $d \in [0,2]$, we have*

$$P(\dim_H \mathcal{H} = d) \in \{0,1\}.$$

**2. The first moment estimate.** Fix $\kappa > 0$ and $z_0 \in \mathbb{H}$; let $\gamma$ be the trace of a chordal $\mathrm{SLE}_\kappa$ in $\mathbb{H}$, and let $\mathcal{H} = \gamma([0, \infty))$ be the image of $\gamma$. We want to compute the probability that $\mathcal{H}$ touches the disk $\mathcal{B}(z_0, \varepsilon)$ for $\varepsilon > 0$.

PROPOSITION 4. *Let $\alpha(z_0) \in (0, \pi)$ be the argument of $z_0$. Then, if $\kappa \in (0, 8)$, we have the following estimate:*

$$P(\mathcal{B}(z_0, \varepsilon) \cap \mathcal{H} \neq \varnothing) \asymp \left( \frac{\varepsilon}{\Im(z_0)} \right)^{1 - \kappa/8} (\sin \alpha(z_0))^{8/\kappa - 1}.$$

*If $\kappa \geq 8$, then this probability is equal to 1 for all $\varepsilon > 0$.*

REMARK. We know that $\mathcal{H}$ is a closed subset of $\bar{\mathbb{H}}$ (indeed, this is a consequence of the transience of $\gamma$—cf. [14]). For $\kappa \geq 8$, this proves that, for all $z \in \bar{\mathbb{H}}$, $P(z \in \mathcal{H}) = 1$, hence, $\mathcal{H}$ almost surely has full measure. And since it is closed, this implies that, with probability 1, $\gamma$ is space-filling, as was already proved by Rohde and Schramm [14] for $\kappa > 8$ and by Lawler, Schramm and Werner [12] for $\kappa = 8$ (for which a separate proof is needed for the existence of $\gamma$).

PROOF OF PROPOSITION 4. The idea of the following proof is originally due to Oded Schramm. Let $\delta_t$ be the Euclidean distance between $z_0$ and $K_t$. $(\delta_t)$ is then a nonincreasing process, and its limit when $t$ goes to $+\infty$ is the distance between $z_0$ and $\mathcal{H}$. Besides, we can apply the Köbe 1/4 theorem to the map $g_t$: this leads to the estimate

$$(2.1) \qquad\qquad \delta_t \asymp \frac{\Im(g_t(z_0))}{|g_t'(z_0)|}$$

(where the implicit constants are universal—namely, 1/4 and 4).

It will be more convenient to fix the image of $z_0$ under the random conformal map. Hence, introduce the following map:

$$\tilde{g}_t : z \mapsto \frac{g_t(z) - g_t(z_0)}{g_t(z) - \overline{g_t(z_0)}}.$$

It is easy to see that $\tilde{g}_t$ maps $\mathbb{H} \setminus K_t$ conformally onto the unit disk $\mathbb{U}$, and maps infinity to 1 and $z_0$ to 0. In other words, the map

$$w \mapsto \tilde{g}_t \left( \frac{w\overline{g_t(z_0)} - g_t(z_0)}{w - 1} \right)$$

maps the complement of some compact $\tilde{K}_t$ in $\mathbb{U}$ onto $\mathbb{U}$, fixing 0 and 1 (in all this proof, $z$ will stand for an element of $\mathbb{H}$ and $w$ for an element of $\mathbb{U}$).

Moreover, in this setup equation (2.1) becomes simpler (because the distance between 0 and the unit circle is fixed): Namely,

$$(2.2) \qquad \delta_t \asymp \frac{1}{|\tilde{g}_t'(z_0)|}.$$

Differentiating $\tilde{g}_t(z)$ with respect to $t$ (which is a little messy and error-prone, but straightforward) leads to the following differential equation:

$$(2.3) \qquad \partial_t \tilde{g}_t(z) = \frac{2(\tilde{\beta}_t - 1)^3}{(g_t(z_0) - \overline{g_t(z_0)})^2 \tilde{\beta}_t^2} \cdot \frac{\tilde{\beta}_t \tilde{g}_t(z)(\tilde{g}_t(z) - 1)}{\tilde{g}_t(z) - \tilde{\beta}_t},$$

where $(\tilde{\beta}_t)$ is the process on the unit circle defined by

$$\tilde{\beta}_t = \frac{\beta_t - g_t(z_0)}{\beta_t - \overline{g_t(z_0)}}.$$

Now the structure of the expression for $\partial_t \tilde{g}_t(z)$ [equation (2.3)] is quite nice: The first factor does not depend on $z$ and the second one only depends on $z_0$ through $\tilde{\beta}$. Hence, let us define a (random) time change by taking the real part of the first factor; namely, let

$$ds = \frac{(\tilde{\beta}_t - 1)^4}{|g_t(z_0) - \overline{g_t(z_0)}|^2 \tilde{\beta}_t^2} \, dt,$$

and introduce $h_s = \tilde{g}_{t(s)}$.

Then equation (2.3) becomes similar to a radial Loewner equation, that is, it can be written as

$$(2.4) \qquad \partial_s h_s(z) = \tilde{X}(\tilde{\beta}_{t(s)}, h_s(z)),$$

where $\tilde{X}$ is the vector field in $\mathbb{U}$ defined as

$$(2.5) \qquad \tilde{X}(\zeta, w) = \frac{2\zeta w(w - 1)}{(1 - \zeta)(w - \zeta)}.$$

The only missing part is now the description of the driving process $\tilde{\beta}$. Applying Itô's formula (now *this* is an ugly computation) and then the previous time-change, we see that $\tilde{\beta}_{t(s)}$ can be written as $\exp(i\alpha_s)$, where $(\alpha_s)$ is a diffusion process on the interval $(0, 2\pi)$ satisfying the equation

$$(2.6) \qquad d\alpha_s = \sqrt{\kappa} \, dB_s + \frac{\kappa - 4}{2} \cot g \frac{\alpha_s}{2} \, ds$$

with the initial condition $\alpha_0 = 2\alpha(z_0)$.

The above construction is licit as long as $z_0$ remains inside the domain of $g_t$. While this holds, differentiating (2.4) with respect to $z$ at $z = z_0$ yields

$$\partial_s h_s'(z_0) = \frac{2h_s'(z_0)}{1 - \tilde{\beta}_s},$$

so that dividing by $h'_s(z_0) \neq 0$ and taking the real parts of both sides we get

$$\partial_s \log |h'_s(z_0)| = 1,$$

that is, almost surely, for all $s > 0$, $|h'_s(z_0)| = |h'_0(z_0)|e^s$. Combining this with (2.2) shows that

$$\delta_{t(s)} \asymp \delta_0 e^{-s} \asymp \Im(z_0)e^{-s}.$$

Finally, let us look at what happens at the stopping time

$$\tau_{z_0} = \mathrm{Inf}\{t : z_0 \in K_t\}.$$

We are in one out of two situations: Either $z_0$ is on the trace: in this case $\delta_t$ goes to 0, meaning that $s$ goes to $\infty$, and the diffusion $(\alpha_s)$ does not touch $\{0, 2\pi\}$. Or, $z_0$ is not on the trace: then $\delta_t$ tends to $d(z_0, \mathcal{H}) > 0$, and the diffusion $(\alpha_s)$ reaches the boundary of the interval $(0, 2\pi)$ at time

$$s_0 := \log \delta_0 - \log d(z_0, \mathcal{H}) + \mathcal{O}(1).$$

Let $S$ be the surviving time of $(\alpha_s)$: the previous construction then shows that

$$d(z_0, \mathcal{H}) \asymp \delta_0 e^{-S},$$

and estimating the probability that $z_0$ is $\varepsilon$-close to the trace becomes equivalent to estimating the probability that $(\alpha_s)$ survives up to time $\log(\delta_0/\varepsilon)$.

Assume for a moment that $\kappa > 4$. The behavior of $\mathrm{cotg}\,\alpha/2$ when $\alpha$ is close to 0 shows that $(\alpha_s)$ can be compared to the diffusion $(\bar{\alpha}_s)$ generated by

$$d\bar{\alpha}_s = \sqrt{\kappa}\,dB_s + (\kappa - 4)\frac{ds}{\bar{\alpha}_s},$$

which (up to a linear time-change) is a Bessel process of dimension

$$d = \frac{3\kappa - 8}{8}.$$

More precisely, $(\bar{\alpha}_s)$ survives almost surely, if and only if $(\alpha_s)$ survives almost surely. But it is known that a Bessel process of dimension $d$ survives almost surely if $d \geq 2$, and dies almost surely if $d < 2$. Hence, we already obtain the phase transition at $\kappa = 8$:

- If $\kappa \geq 8$, then $d \geq 2$, and $(\alpha_s)$ survives almost surely. Hence, almost surely, $d(z_0, \mathcal{H}) = 0$, and for all $\varepsilon > 0$, the trace will almost surely touch $\mathcal{B}(z_0, \varepsilon)$.
- If $\kappa < 8$, then $d < 2$ and $(\alpha_s)$ dies almost surely in finite time. Hence, almost surely, $d(z_0, \mathcal{H}) > 0$.

So, there is nothing left to prove for $\kappa \geq 8$. From now on, we shall then suppose that $\kappa \in (0, 8)$. If $\kappa \leq 4$, then the drift of $(\alpha_s)$ is toward the boundary, hence, comparing it to standard Brownian motion shows that it dies almost surely in finite time as for $\kappa \in (4, 8)$. We want to apply Lemma 2 to $(\alpha_s)$ and for that we need to know the principal eigenvalue of the generator $L_\kappa$ of the diffusion. It can be seen that the function

$$(\sin(x/2))^{8/\kappa - 1}$$

is a positive eigenfunction of $L_\kappa$, with eigenvalue $1 - \kappa/8$: hence, we already obtain that, if $\alpha_0$ is far from the boundary, $P(S > s) \asymp \exp(-(1 - \kappa/8)s)$, that is,

$$(2.7) \qquad P(d(z_0, \mathcal{H}) \leq \varepsilon) \asymp e^{(1 - \kappa/8) \log(\varepsilon/\delta_0)} \asymp \left( \frac{\varepsilon}{\delta_0} \right)^{1 - \kappa/8},$$

which is the correct estimate. It remains to take the value of $\alpha_0$ into account.

Introduce the following process:

$$X_s := \sin\left( \frac{\alpha_s}{2} \right)^{8/\kappa - 1} e^{(1 - \kappa/8)s}$$

(and $X_s = 0$ if $s \geq S$). Applying the Itô formula shows that $(X_s)$ is a local martingale [in fact, this is the same statement as saying that $\sin(x/2)^{8/\kappa - 1}$ is an eigenfunction of the generator], and it is bounded on any bounded time interval. Hence, taking the expected value of $X$ at times $0$ and $s$ shows that

$$(2.8) \qquad \sin\left( \frac{\alpha_0}{2} \right)^{8/\kappa - 1} = e^{(1 - \kappa/8)s} P(S \geq s) E\left[ \sin\left( \frac{\alpha_s}{2} \right)^{8/\kappa - 1} \Big| S \geq s \right].$$

The same proof as that of Lemma 2 shows that, for all $s \geq 1$,

$$P(\alpha_s \in [\pi/2, 3\pi/2] | S \geq s) > 0,$$

with constants depending only on $\kappa$; combining this with $(2.8)$ then provides

$$P(S \geq s) \asymp e^{-(1 - \kappa/8)s} \sin\left( \frac{\alpha_0}{2} \right)^{8/\kappa - 1},$$

again with constants depending only on $\kappa$. Applying the same computation as for equation $(2.7)$ ends the proof. $\square$

COROLLARY 5. *Let $D \subsetneq \mathbb{C}$ be a simply connected domain, $a$ and $b$ be two points on the boundary of $D$, and $\gamma$ be the path of a chordal $\mathrm{SLE}_\kappa$ in $D$ from $a$ to $b$, with $\kappa \in (0, 8)$. Then, for all $z \in D$ and $\varepsilon < d(z, \partial D)/2$, we have*

$$P(\gamma \cap \mathcal{B}(z, \varepsilon) \neq \varnothing) \asymp \left( \frac{\varepsilon}{d(z, \partial D)} \right)^{1 - \kappa/8} (\omega_z(ab) \wedge \omega_z(ba))^{8/\kappa - 1},$$

*where $\omega_z$ is the harmonic measure on $\partial D$ seen from $z$ and $ab$ is the positively oriented arc from $a$ to $b$ along $\partial D$.*

Proof. This is easily seen by considering a conformal map $\Phi$ mapping $D$ to the upper-half plane, $a$ to 0 and $b$ to $\infty$: Since the harmonic measure from $z$ in $D$ is mapped to the harmonic measure from $\Phi(z)$ in $\mathbb{H}$, it is sufficient to prove that, for all $z \in \mathbb{H}$,

$$\omega_z(\mathbb{R}_+) \wedge \omega_z(\mathbb{R}_-) \asymp \sin(\arg z);$$

and $\omega_z(\mathbb{R}_+)$ can be explicitly computed, because $\omega_z$ is a Cauchy distribution on the real line:

$$\omega_{x+iy}(\mathbb{R}_+) = \frac{1}{\pi} \int_0^\infty \frac{du/y}{1+(u-x)^2/y^2} = \frac{1}{2} + \frac{1}{\pi} \operatorname{arctg}(x/y).$$

When $x$ tends to $-\infty$, this behaves like $-y/\pi x$, which is equivalent to $\sin(\arg(x+iy))/\pi$. □

This "intrinsic" formulation of the hitting probability will make the derivation of the second moment estimate more readable.

**3. The second moment estimate.** We still have to derive condition 3 in Proposition 1. For $\kappa = 6$, it was obtained using the locality property, but this does not hold for other values of $\kappa$, so we can rely only on the Markov property. In this whole section we shall assume that $\kappa < 8$ (there is nothing to prove if $\kappa \geq 8$, since in that case $\gamma$ is space-filling).

The general idea is as follows. Fix two points $z$ and $z'$ in the upper half plane, and $\varepsilon < |z' - z|/2$. We want to estimate the probability that the trace $\gamma$ visits both $\mathcal{B}(z, \varepsilon)$ and $\mathcal{B}(z', \varepsilon)$. Assume that it touches, say, the first one (which happens with probability of order $\varepsilon^{1-\kappa/8}$), and that it does so before touching the other.

Apply the Markov property at the first hitting time $T_\varepsilon(z)$ of $\mathcal{B}(z, \varepsilon)$: If everything is going fine and we are lucky, the distance between $z'$ and $K_{T_\varepsilon(z)}$ will still be of order $|z' - z|$. Hence, applying the first moment estimate to this situation shows that the conditional probability that $\gamma$ hits $\mathcal{B}(z', \varepsilon)$ is not greater than $C(\varepsilon/|z' - z|)^{1-\kappa/8}$ [it might actually be much smaller, if the real part of $g_{T_\varepsilon(z)}(z')$ is large, but this is not a problem since we only need an upper bound], and this gives the right estimate for the second moments:

$$C \frac{\varepsilon^{2-\kappa/4}}{|z' - z|^{1-\kappa/8}}.$$

The whole point is then to prove that this is the main contribution to the second moment; the way we achieve it is by providing sufficiently sharp upper bounds for the second term of the estimate given by Corollary 5.

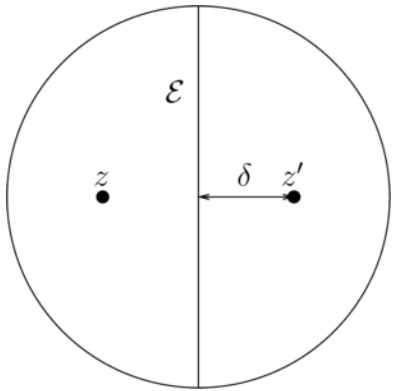

Fig. 1. *Second moments: the setup.*

3.1. *Preliminaries.* The first part of the proof is a succession of topological lemmas which allow for a precise estimation of the harmonic measures of the two sides of the *SLE* process. They are easier to state in the case $\kappa \leq 4$, for which the process consists in a simple curve. In the case $4 < \kappa < 8$, what happens is that a positive area is "swallowed" by the process; in all the following discussion, nothing changes as long as the points $z$ and $z'$ themselves are not swallowed, and the arguments are exactly the same—as all that is required for the proofs to apply is for the complement of the process to be simply connected and contain both $z$ and $z'$.

On the other hand, if (say) $z$ is swallowed at a given time, at which the curve has not touched $\mathcal{C}(z, \varepsilon)$ yet, then this will never happen, so this event does not contribute at all to the probability of the event we are interested in. If the trace does touch $\mathcal{C}(z, \varepsilon)$ before swallowing $z$, then the swallowing occurs at a time when it is not relevant anymore—since we know already that $z$ is $\varepsilon$-close to the path—so again the rest of the argument is not affected.

In order to simplify the exposition of the argument, we will implicitly assume that indeed $\kappa \leq 4$. The interested reader can easily check as she proceeds that what follows does apply to the other cases, with little change in the writing.

Let $z$, $z'$ be two points in the upper half plane, and let $\delta = |z - z'|/2$. We can assume that both $\Im z$ and $\Im z'$ are greater than $18\delta$ (say). Introduce a "separator set" (cf. Figure 1):

$$\mathcal{E} = \mathcal{C}\left(\frac{z + z'}{2}, 2\delta\right) \cup \{w \in \mathbb{H} : d(w, z) = d(w, z') \leq \delta\sqrt{5}\}.$$

At each positive time $t$, the complement $H_t$ of $K_t$ in $\mathbb{H}$ is an open and simply connected domain, hence, its intersection with $\mathcal{E}$ is the disjoint union of at most countably many connected sets, each separating $H_t$ into two (or up to four for at most two of them) connected components. If both $z$ and $z'$ are in $H_t$, let $\mathcal{E}_t$ be the union of those crosscuts which disconnect $z$ from $z'$ in $H_t$; if either $z$ or $z'$ is in $K_t$, let $\mathcal{E}_t = \varnothing$—notice that in the case of an SLE process with parameter $\kappa \leq 4$, this almost surely never happens. Note that, as long as $z$ and $z'$ are in $H_t$, $\mathcal{E}_t$ is not empty, because $H_t$ is simply connected and $\mathcal{E}$ itself disconnects $z$ from $z'$. The components of $\mathcal{E}_t$ can then be ordered in the way they first appear on any path going from $z$ to $z'$ in $H_t$; let $\lambda_t$ be the first one, and $\lambda'_t$ the last one (which is also the first one seen from $z'$ to $z$); for convenience, in the case $\mathcal{E}_t = \varnothing$, let $\lambda_t = \lambda'_t = \varnothing$ too.

For each time $t$ (possibly random), introduce

$$(3.1) \qquad \tilde{t} := \mathrm{Inf}\{s > t : K_s \cap \lambda'_t \neq \varnothing\},$$

$$(3.2) \qquad \check{t} := \mathrm{Inf}\{s > t : K_s \cap \lambda_t \neq \varnothing\},$$

with the usual convention that the infimum of the empty set is infinite. Clearly, $\tilde{\tau}$ and $\check{\tau}$ are stopping times if $\tau$ is one.

Besides, $\mathcal{E}_t$ does not change on any time-interval on which $\gamma$ does not intersect $\mathcal{E}$—hence, if for some $t_1 < t_2$, $\gamma((t_1, t_2)) \cap \mathcal{E} = \varnothing$, we can define $\mathcal{E}_{t_2-}$ as its constant value on the interval $(t_1, t_2)$ (i.e., as $\mathcal{E}_{t_1}$). We say that a positive stopping time $t$ is a *good time* if the following conditions are satisfied with probability 1:

- There exists $s < t$ such that $\gamma((s, t)) \cap \mathcal{E} = \varnothing$;
- $\gamma(t) \in \mathcal{E}_{t-}$.

Examples of good times are $\tilde{t}$ and $\check{t}$ when $t$ is a stopping time such that $\gamma(t) \notin \mathcal{E}$ holds with probability 1.

We first give two preliminary lemmas which will be useful in estimating the harmonic measures appearing in the statement of Corollary 5. They are not specific to SLE, but they depend (as stated) on the fact that $\gamma$ is a simple curve which does not contain $z$ nor $z'$ (as is the case with probability 1 in the case $\kappa \leq 4$); they have obvious counterparts obtained by exchanging $z$ and $z'$ and replacing everywhere $\tilde{t}$ with $\check{t}$.

For each positive time, let $\omega_t$ (resp. $\omega'_t$) be the smaller of the harmonic measures of the two sides of $\gamma$, from $z$ (resp. $z'$) in $H_t$—this corresponds to the term we want to estimate in the statement of Corollary 5. (Here and in all the sequel, as is natural, we include the positive real axis in the right side of $\gamma$ and the negative real axis in its left side.) Besides, for $t > 0$ and $\rho \in (0, \delta)$, let $\mathcal{B}_t(\rho)$ [resp. $\mathcal{B}'_t(\rho)$] be the closure of the connected component of $z$ (resp. $z'$) in $\mathcal{B}(z, \rho) \cap H_t$ [resp. $\mathcal{B}(z', \rho) \cap H_t$].

In all that follows we will use the following notation at each time $t > 0$ (together with their counterparts around $z'$):

$$r_t := d(z, \gamma([0, t]) \cup \mathbb{R});$$

$$\rho_t := \mathrm{Inf}\{\rho \in (0, \delta) : \mathcal{B}_t(\rho) \text{ disconnects } z' \text{ from } \infty \text{ in } H_t\}$$

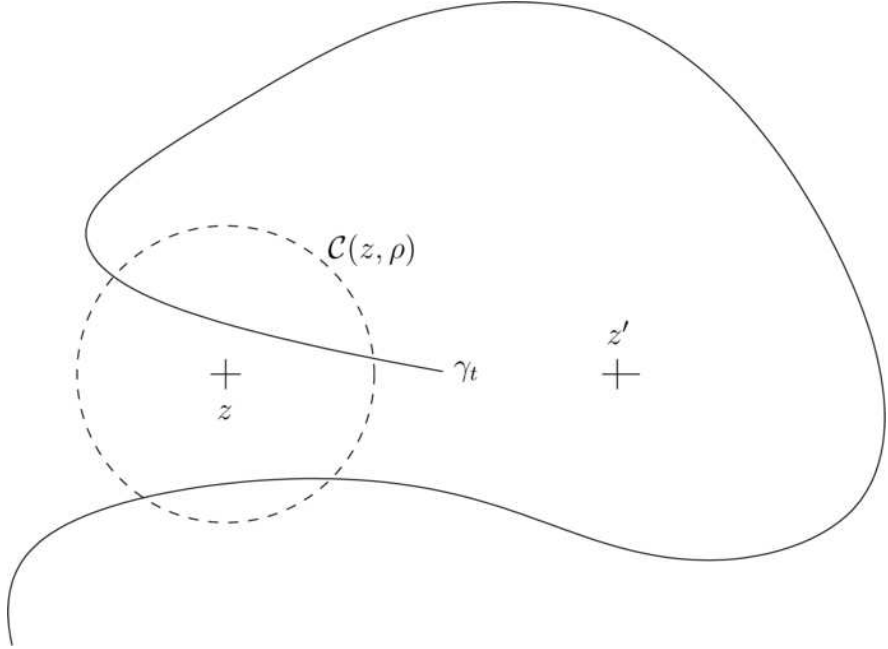

Fig. 2. *Proof of Lemma 6.*

(letting $\rho_t = \delta$ if the infimum is taken over an empty set). Obviously $(r_t)$ is nonincreasing; but $(\rho_t)$ is not in general. Besides, $\rho_t \geq r_t$. Last, since the sets $\mathcal{B}_t(\rho)$ and $\gamma([0, t])$ are all compact, it is easy to see that at each time $t$ such that $\rho_t < \delta$, $\gamma([0, t]) \cup \mathcal{B}_t(\rho_t)$ itself does disconnect $z'$ from infinity.

LEMMA 6. *There exists a positive constant $c$ such that the following happens. Let $t$ be a good time, and $\rho \in (r_t, \delta)$. If $\omega_t \geq c(r_t/\rho)^{1/2}$, then $\gamma([0,t]) \cup \mathcal{B}_t(\rho)$ disconnects $z'$ from infinity; in particular, $\rho_t \leq \rho$.*

PROOF. First make the following remark: For $\rho \in (r_t, \delta)$, if the harmonic measure from $z$ in $\mathcal{B}_t(\rho) \setminus \gamma([0, t])$ gives positive mass to both sides of $\gamma$, then $\mathcal{B}_t(\rho)$ separates $z'$ from infinity in $H_t$ (see Figure 2).

Indeed, assume that we are this case. That means that there exist two disjoint smooth curves $\zeta_1, \zeta_2 : [0, 1] \to \mathbb{H}$ satisfying $\zeta_1(0) = \zeta_2(0) = z$, $\zeta_i((0, 1)) \subset \mathcal{B}_t(\rho) \setminus \gamma([0, t])$ and $\zeta_i(1) \in \gamma([0, t])$, each landing on a different side of $\gamma$ [i.e., $\lim_{s \to 1} g_t(\zeta_1(s)) - \beta_t \in (0, +\infty)$ and $\lim_{s \to 1} g_t(\zeta_2(s)) - \beta_t \in (-\infty, 0)$—note that such limits are always welldefined because $g_t$ extends continuously to the boundary of $H_t$]. Let $\zeta = \zeta_1((0, 1)) \cup \zeta_2((0, 1)) \cup \{z\}$ be the corresponding cross-cut: The complement of $\zeta$ in $H_t$ has exactly two (simply) connected components, one of which is unbounded.

If $\gamma(t)$ were on the boundary of the unbounded component, then one could continue $\gamma([0, t])$ with some curve $\hat{\gamma}$ contained in $H_t \setminus \zeta$ and tending to infinity. Then the bounded component of $H_t \setminus \zeta$ would be contained in one of the components of $\mathbb{H} \setminus (\gamma([0, t]) \cup \hat{\gamma})$, hence, its boundary (which contains both endpoints of $\zeta$) would intersect only one side of $\gamma$—which is in contradiction with our hypothesis.

Now if $z'$ were in the unbounded component, it would be possible to join $z$ to $z'$ inside the unbounded component. But such a path would have to intersect the part of $\mathcal{E}_{t-}$ which contains $\gamma(t)$ (by the definition of $\mathcal{E}_{t-}$ and that of a good time), and it does not because this part of $\mathcal{E}_t$ is contained in the connected component of $H_t \setminus \zeta$ which contains $\gamma(t)$ on its boundary—that is, the bounded one.

To sum it up, $\gamma([0, t]) \cup \zeta$ cuts $H_t$ into two connected components, and the bounded component contains $z'$ in its interior and $\gamma(t)$ on its boundary. In particular, since $\zeta$ is contained in $\mathcal{B}_t(\rho)$, this implies that $\mathcal{B}_t(\rho)$ separates $z'$ from infinity in $H_t$.

It is then straightforward to complete the proof of the lemma, by applying Beurling's estimate in $\mathcal{B}_t(\rho)$ and the maximum principle. $\quad\square$

We will actually use the converse of this lemma: At any good time $t$, we have

$$(3.3) \qquad \omega_t \leq c\left(\frac{r_t}{\rho_t}\right)^{1/2}.$$

A related fact is the following:

LEMMA 7. Let $r \in (0, \delta)$, $T$ be the first time when $\gamma$ hits the circle $\mathcal{C}(z, r)$ and $\tilde{T}$ as introduced above [see equation (3.1)]. If $\tilde{T}$ if finite, then $\mathcal{B}_{\tilde{T}}(r)$ does not disconnect $z'$ from infinity in $H_{\tilde{T}}$.

PROOF. Let $\zeta$ be a continuous, simple curve going from $z'$ to infinity in $H_{\tilde{T}}$, and let $\mathcal{C}$ be the connected component of $H_T \setminus \mathcal{E}$ which contains $z'$. It is always possible to ensure that $\zeta$ intersect every component of $\mathcal{E}_{\tilde{T}}$ at most once. The boundary of $\mathcal{C}$ is contained in $\gamma([0, T]) \cup \lambda'_T$; and since $\gamma([0, T]) \cap \mathcal{E}$ is not empty [because $\mathcal{E}$ separates $\mathcal{C}(z, r)$ from 0], necessarily $\gamma([0, T]) \cap \bar{\lambda}'_T \neq \varnothing$. In particular, $H_T \setminus (\zeta \cup \mathcal{C} \cup \lambda'_T)$ has exactly two unbounded connected components, say, $U_1$ and $U_2$.

Assume that $\mathcal{B}_{\tilde{T}}(r)$ does disconnect $z'$ from infinity. Then $\zeta$ has to intersect its interior, splitting it into at most countably many connected components. $\gamma([0, \tilde{T}])$ has to intersect at least one component on each side of $\zeta$, since if not one could deform $\zeta$ so that it avoids $\mathcal{B}_{\tilde{T}}(r)$—but by the definition of $T$, $\gamma([0, T])$ intersects only the adherence of one component, say, on the left of $\zeta$. By construction, the only way for $\gamma$ to reach a component on the other side of $\zeta$ is by intersecting $\mathcal{C}$ and hence $\lambda'_T$, so it cannot happen before time $\tilde{T}$. $\quad\square$

In other words, if $T$ is the first time when $\gamma$ intersects $\mathcal{C}(z, r)$ and $\tau$ is the first time $t$ such that $\rho_t \leq r$, then assuming that $\tau$ is finite, we have $T < \tilde{T} < \tau$.

The last lemma in this section is specific to SLE: It is a quantitative version of the transience of the curve $\gamma$ and basically says that if $\gamma$ forms a fjord, then it is not likely to enter it. With the modifications of notation described later for the case $4 < \kappa < 8$, it holds also in that case, and the proof is the same.

Lemma 8. *Let $\gamma$ be the trace of an SLE with parameter $\kappa \leq 4$; then there exist positive constants $C$ and $\eta$ such that the following happens. Let $\rho > 0$ and let $\tau$ be the first time $t$ such that $\rho_t \leq \rho$ (i.e., such that $\gamma([0, t]) \cup \mathcal{B}_t(\rho)$ disconnects $z'$ from infinity). $\tau$ is finite with positive probability, in which case we have $|\gamma_\tau - z| = \rho$, and:*

1. *$P(\tilde{\tau} < \infty | \mathcal{F}_\tau, \tau < \infty) \leq C(\rho/\delta)^\eta$;*
2. *For every $r < r_\tau$,*

   *$P(\tilde{\tau} < \infty, r_{\tilde{\tau}} < r | \mathcal{F}_\tau, \tau < \infty) \leq C(r/r_\tau)^{1-\kappa/8}(\rho/\delta)^\eta.$*

Proof. (i) In all this proof $c$ will denote any finite positive constant which depends only on $\kappa$. Notice that if $z' \in K_\tau$ (which can happen if $\kappa > 4$), there is nothing to prove, since $\tilde{\tau} = \infty$ in this case; so we assume from now on that $z' \notin K_\tau$. Recall that $\lambda'_\tau$ is the last component of $\mathcal{E}_\tau$ that one has to cross when going from $z$ to $z'$ in $H_\tau$. By monotonicity, the extremal distance between $\mathcal{E}$ and $\mathcal{B}_\tau(\rho)$ in $H_\tau$ is bounded below by $\frac{1}{2\pi}\log(\delta/\rho)$—and hence so is the extremal distance between $\lambda'_\tau$ and $\mathcal{B}_\tau(\rho)$.

By the definition of $\tau$, it is possible to find a simple continuous curve $\zeta$ going from $z'$ to $\infty$ in $(H_\tau \setminus \mathcal{B}_\tau(\rho)) \cup \{\gamma_\tau\}$ (e.g., by using the fact that $\gamma([0, \tau - s]) \cup \mathcal{B}_\tau(\rho)$ does *not* separate $z'$ from $\infty$ for $s > 0$, choosing $\zeta^s$ accordingly, and letting $s$ go to 0); and there exists a simple continuous curve $\zeta'$ going from $z$ to $\gamma_\tau$ in $\mathcal{B}_\tau(\rho)$ (see Figure 3). Considering these curves, it is easy to see that $\lambda'_\tau$ disconnects $z'$, and not $z$, from infinity in $H_\tau$.

The construction also shows that $|\gamma_\tau - z| \leq \rho$: Indeed, if not then, we can deform $\zeta$ locally around $\gamma_\tau$ to obtain a continuous curve going from $z'$ to infinity without hitting $\gamma([0, \tau]) \cup \mathcal{B}_\tau(\rho)$, which is in contradiction with the definition of $\tau$. If on the other hand we had $|\gamma_\tau - z| < \rho$, then for $s < \tau$ large enough, $\mathcal{B}_s(\rho) \cup \gamma([0, s])$ would still disconnect $z'$ from infinity, also leading to a contradiction—hence, $|\gamma_\tau - z| = \rho$.

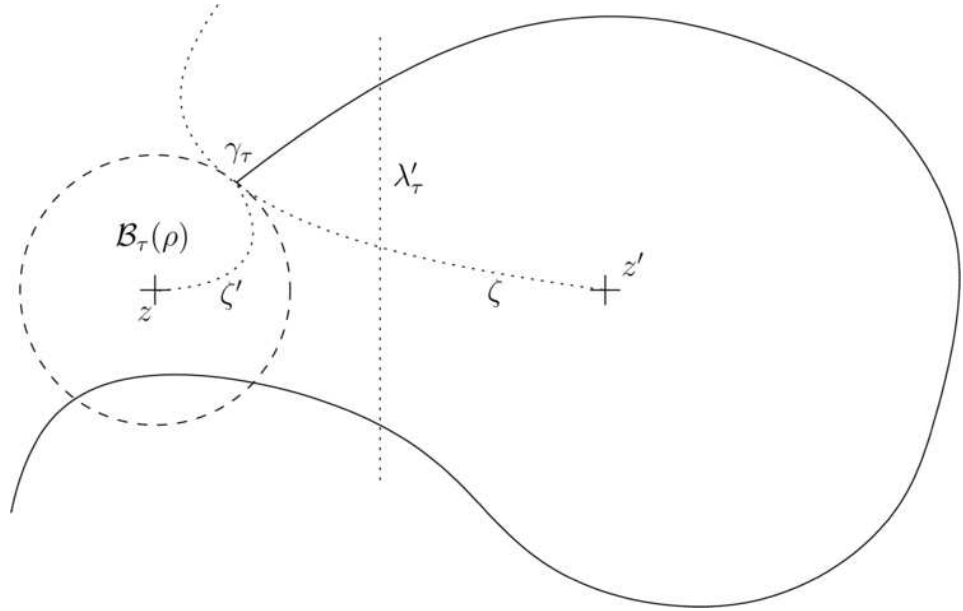

FIG. 3.   *Proof of Lemma* 8: *Setup.*

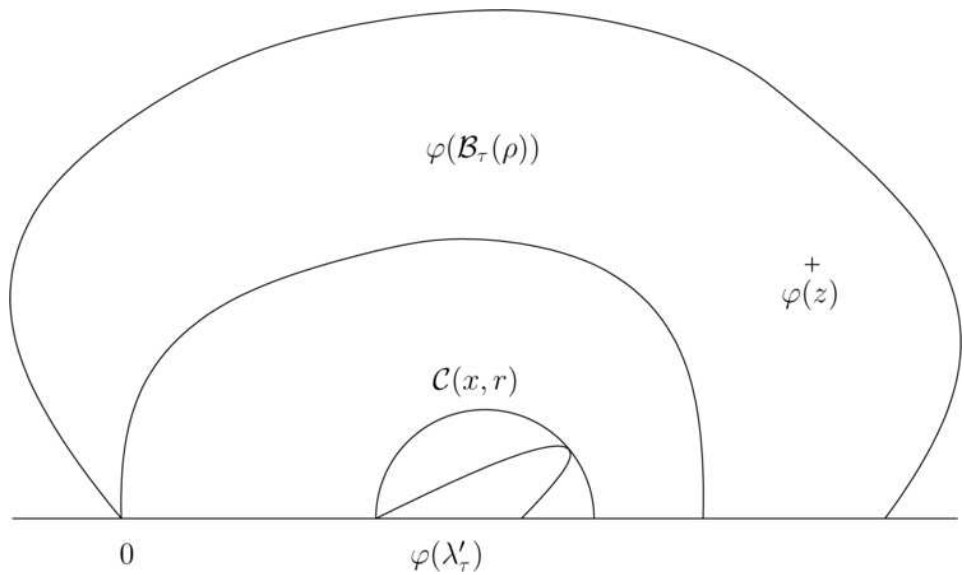

FIG. 4.   *Proof of Lemma* 8: *After mapping by* $\varphi$.

Map the whole picture by $\varphi := g_\tau - \beta_\tau$. $\mathcal{B}_\tau(\rho)$ is mapped to a cross-cut having 0 on its closure, and the images of $\lambda'_\tau$ and $z'$ are in a bounded connected component of its complement (cf. Figure 4). The boundary of the unbounded component of the complement of $\varphi(\mathcal{B}_\tau(\rho))$ contains either $(-\infty, 0)$ or $(0, \infty)$; for ease of notation, we assume that the former holds, as in the figure.

By conformal invariance, the extremal distance in $\mathbb{H}$ between $\varphi(\lambda'_\tau)$ and $\varphi(\mathcal{B}_\tau(\rho))$ is bounded below by $\frac{1}{2\pi} \log(\delta/\rho)$—and so is the extremal distance between $\varphi(\lambda'_\tau)$ and $(-\infty, 0)$, since 0 is on the closure of $\varphi(\mathcal{B}(z, \rho))$. Let $x > 0$ be the smallest element of $\mathbb{R} \cap \overline{\varphi(\lambda'_\tau)}$ and let $r > 0$ be the smallest positive real such that $\mathcal{C}(x, r)$ separates $\varphi(\lambda'_\tau)$ from infinity [so that, in particular, $\varphi(\lambda'_\tau) \subset \bar{\mathcal{B}}(x, r)$]. Let $l$ be the extremal distance in $\mathbb{H} \setminus \varphi(\lambda'_\tau)$ between $(0, x)$ and $(\mathrm{Max}(\mathbb{R} \cap \overline{\varphi(\lambda'_\tau)}), +\infty)$, that is, the reciprocal of the extremal distance in $\mathbb{H}$ between $(-\infty, 0)$ and $\varphi(\lambda'_\tau)$: We have $l \leq c/\log(\delta/\rho)$.

On the other hand, it is possible to find a lower bound for $l$ in terms of $x$ and $r$, as follows: Consider the metric given by $u(z) = \alpha/|z - x|$ if $r/5 < |z - x| < 5x$, and $u(z) = 0$ otherwise, where $\alpha$ is chosen to normalize the surface integral to 1—so that $\alpha$ is of order $1/\log(x/r)$. In this metric, the length of any curve joining $(0, x)$ to $(\mathrm{Max}(\mathbb{R} \cap \overline{\varphi(\lambda'_\tau)}), +\infty)$ around $\varphi(\lambda'_\tau)$ is at least of order $\alpha$ [this can be seen, e.g., by using the conformal map $z \mapsto \log(z - x)$, which maps $u$, where it is not zero, to the renormalized Euclidean metric]. Hence, we obtain a lower bound for $l$, of the form $l \geq c/\log(x/r)$, and thus an upper bound on $r/x$ of the form $c(\rho/\delta)^\alpha$ for some $\alpha > 0$.

Let $p(r)$ be the probability that a chordal $\mathrm{SLE}_\kappa$ starting at 0 touches the circle $\mathcal{C}(1, r)$. Since we are in the case $\kappa < 8$, $0 < p(r) < 1$ as soon as $r \in (0, 1)$, and $p(r)$ goes to 0 with $r$; besides, the strong Markov property applied at the first hitting time of $\mathcal{C}(1, r)$ (if such a time exists) together with Köbe's 1/4-theorem ensure that there is a finite positive constant $C > 1$ such that, for all $r, r' < 1$, $p(rr') \leq Cp(r)p(r')$. So, let $r_0$ be such that $p(r_0) < C^{-2}$ and apply the inequality $n - 1$ times to obtain $p(r_0^n) \leq C^{n-1}(C^{-2})^n \leq C^{-n}$. This implies that $p(r)$ is bounded above by $cr^{\eta'}$ for some $\eta' > 0$ [actually, the optimal value for $\eta'$ is the same as the boundary exponent $s_b = (8/\kappa) - 1$, but we will not need this].

Hence, we obtain

$$P(\tilde{\tau} < \infty | \mathcal{F}_\tau, \tau < \infty) \leq c \left(\frac{r}{x}\right)^{\eta'} \leq c \left(\frac{\rho}{\delta}\right)^\eta,$$

with $\eta > 0$, as we wanted.

(ii) The proof of this estimate is actually a simpler version of the proof of the second-moment estimate in the next section, so we will explain it in more detail than would probably be necessary. We want to estimate the conditional probability, conditionally to $\mathcal{F}_\tau$, that $\tilde{\tau}$ is finite and that the curve $\gamma$ hits the circle $\mathcal{C}(z, r)$ before $\tilde{\tau}$ (we say that $\gamma$ *succeeds* if these two conditions are satisfied). Fix $a \in (0, 1)$ (its value will be chosen later in the proof): If $\gamma$ succeeds, then, in particular, it has to hit all the circles of the form $\mathcal{C}(z, r_\tau a^k)$ lying between $\gamma(\tau)$ and $\mathcal{C}(z, r)$, and (the relevant parts of) all the circles of the form $\mathcal{C}(z, \rho a^{-k})$ lying between $\gamma(\tau)$ and $\mathcal{E}_\tau$.

The idea is then the following: For each possible ordering of these hitting times, we will estimate the probability that the circles are hit in this particular order, using the strong Markov property recursively together with previous estimates; we can then sum over all possible orderings to obtain an estimate of the probability that $\gamma$ succeeds.

For each $k > 0$, let $T_k$ be the first hitting time of $\mathcal{C}(z, r_\tau a^k)$ by $\gamma$. Besides, let $\lambda_k$ be the last connected component of $\mathcal{C}(z, \rho a^{-k}) \cap H_\tau$ which a curve going from $z$ to $z'$ has to cross, and let $\tau_k$ be the first hitting time of $\lambda_k$ by $\gamma$. Last, let $k_1$ (resp. $k_2$) be the largest integer smaller than $\log(\delta/\rho)/\log(1/a)$ [resp. $\log(r_\tau/r)/\log(1/a)$]: It is sufficient to give an upper bound for the probability that both $\tau_{k_1}$ and $T_{k_2}$ are finite and smaller that $\tilde{\tau}$.

We describe the ordering of the hitting times by specifying the successive numbers of circles of each kind which $\gamma$ hits before time $\tilde{\tau}$. More precisely, assume that $\gamma$ succeeds: Then there are nonnegative integers $I$, $(m_i)_{i \leq I}$ and $(l_i)_{i \leq I}$, all positive except possibly for $m_1$ and $l_I$, such that

$$\tau_1 < \cdots < \tau_{m_1} < T_1 < \cdots < T_{l_1} < \tau_{m_1+1} < \cdots < \tau_{m_1+m_2} < \cdots < T_{l_1+\cdots+l_I} < \tilde{\tau}$$

and $\sum m_i = k_1$ and $\sum l_i = k_2$ [so that $\gamma$ first crosses $m_1$ of the $\lambda_k$, then $l_1$ of the $\mathcal{C}(z, r_\tau a^k)$, then $m_2$ new $\lambda_k$, etc.].

Notice that at time $\tau_i$, Beurling's estimate in the domain $\mathcal{B}(z, r_\tau) \setminus \gamma([0, \tau_i])$ shows that $\omega_{\tau_i}$ is at most equal to $C(r_{\tau_i}/r_\tau)^{1/2}$ (by Lemma 7 and the same argument as the one used in the proof of Lemma 6). Besides, the same proof as that of point (i) in the present lemma shows the following: For given values of the $(m_i)$ and $(l_i)$, for each $i$, conditionally to $\mathcal{F}_{T_{l_1+\cdots+l_i}}$ and the facts that $T_{l_1+\cdots+l_i} < \infty$ and that the last $\tau$-time happening before $T_{l_1+\cdots+l_i}$ is $\tau_{m_1+\cdots+m_i}$ (as will be the case in the construction), $P(\tau_{m_1+\cdots+m_{i+1}} < \infty | \mathcal{F}_{T_{l_1+\cdots+l_i}})$ is bounded by $Ca^{\eta m_{i+1}}$.

For given values of $I$ and the $(m_i)$ and $(l_i)$, applying the strong Markov property at each of the times $T_{l_1+\cdots+l_i}$ and $\tau_{m_1+\cdots+m_i}$, we get an estimate of the conditional probability (conditionally to $\mathcal{F}_\tau$) that $\gamma$ succeeds with this particular ordering, as a product of conditional probabilities, namely,

$$\prod_{i=1}^{I} P(\tau_{m_1+\cdots+m_i} < \infty | \mathcal{F}_{T_{l_1+\cdots+l_{i-1}}}) P(T_{l_1+\cdots+l_i} < \infty | \mathcal{F}_{\tau_{m_1+\cdots+m_i}}).$$

Using the previous estimates, and Corollary 5, this product is bounded above by

$$\prod_{i=1}^{I} Ca^{\eta m_i}(a^{l_1+\cdots+l_{i-1}})^{(8/\kappa-1)/2}(a^{l_i})^{1-\kappa/8}.$$

It remains to sum this estimate over all possible values of $I$, the $m_i$ and the $l_i$. We get the following, where as is usual $C$ is allowed to change from line to line, but depends only on $\kappa$ and later on $a$:

$$P(\tau_{k_1} < \tilde{\tau}, T_{k_2} < \tilde{\tau}) \leq \sum_{I=1}^{\infty} \sum_{(m_i),(l_i)} \prod_{i=1}^{I} C a^{\eta m_i} (a^{l_1 + \cdots + l_{i-1}})^{(8/\kappa-1)/2} (a^{l_i})^{1-\kappa/8}$$

$$\leq a^{\eta k_1 + (1-\kappa/8)k_2} \sum_{I=1}^{\infty} C^I \sum_{(m_i),(l_i)} \prod_{i=1}^{I} (a^{l_1 + \cdots + l_{i-1}})^{(8/\kappa-1)/2}$$

$$= a^{\eta k_1 + (1-\kappa/8)k_2} \sum_{I=1}^{\infty} C^I \sum_{(m_i),(l_i)} \prod_{i=1}^{I} (a^{(I-i)l_i})^{(8/\kappa-1)/2}.$$

For a fixed value of $I$, the number of possible choices for the $m_i$ (which are $I$ integers of sum $k_1$) is smaller than $2^{I+k_1}$, hence, replacing $C$ by $2C$, we get

$$P(\tau_{k_1} < \tilde{\tau}, T_{k_2} < \tilde{\tau}) \leq a^{\eta k_1 + (1-\kappa/8)k_2} 2^{k_1} \sum_{I=1}^{\infty} C^I \sum_{(l_i)} \prod_{i=1}^{I} (a^{(I-i)l_i})^{(8/\kappa-1)/2}.$$

The sum over $(l_i)$ is taken over all $I$-tuples of positive integers with sum $k_2$, so if the first $I-1$ are known, so is the last one. An upper bound is then given by relaxing the condition $l_1 + \cdots + l_I = k_2$ and simply summing over all positive values of $l_1, \ldots, l_{I-1}$ ($l_I$ does not contribute to the product anyway). So we obtain

$$P(\tau_{k_1} < \tilde{\tau}, T_{k_2} < \tilde{\tau}) \leq a^{\eta k_1 + (1-\kappa/8)k_2} 2^{k_1} \sum_{I=1}^{\infty} C^I \prod_{i=1}^{I-1} \sum_{l>0} (a^{(I-i)l})^{(8/\kappa-1)/2}.$$

We can sum over $l > 0$ in each term of the product; each sum will be equal to $a^{(I-i)(8/\kappa-1)/2}$ up to a constant which, if $a$ is chosen small enough, is smaller than 2. Hence, a factor $2^I$ which can again be made part of $C^I$ by doubling the value of $C$:

$$P(\tau_{k_1} < \tilde{\tau}, T_{k_2} < \tilde{\tau}) \leq a^{\eta k_1 + (1-\kappa/8)k_2} 2^{k_1} \sum_{I=1}^{\infty} C^I \prod_{i=1}^{I-1} (a^{(I-i)})^{(8/\kappa-1)/2}.$$

Compute the product explicitly: The exponent of $a$ is then the sum of the $I - i$ for $1 \leq i \leq I - 1$, which is equal to $I(I-1)/2$. The linear term $-I/2$ can be incorporated in the factor $C^I$ (making $C$ depend on $a$ now, which will not be a problem), leading to

$$P(\tau_{k_1} < \tilde{\tau}, T_{k_2} < \tilde{\tau}) \leq a^{\eta k_1 + (1-\kappa/8)k_2} 2^{k_1} \sum_{I=1}^{\infty} C^I (a^{I^2/2})^{(8/\kappa-1)/2}$$

$$\leq (2a^{\eta/2})^{k_1} a^{(\eta/2)k_1 + (1 - \kappa/8)k_2} \sum_{I=1}^{\infty} C^I (a^{I^2/2})^{(8/\kappa - 1)/2}.$$

Now pick $a$ small enough that $2a^{\eta/2}$ is smaller than 1. The sum in the previous expression is finite (because $\kappa < 8$ and $a < 1$), so we obtain

$$P(\tau_{k_1} < \tilde{\tau}, T_{k_2} < \tilde{\tau}) \leq Ca^{(\eta/2)k_1 + (1 - \kappa/8)k_2},$$

which implies the announced result. $\square$

REMARK 1.   It is possible to simplify the statement of the last part of the proof of the lemma (though unfortunately not the computation) in the following way. Let $\mathbf{m} = (m_i)$ and $\mathbf{l} = (l_i)$ be the jump sizes of the process, which we will interpret as ordered partitions of $k_1$ and $k_2$, respectively. As is customary, we write this as $\mathbf{m} \vdash k_1$, respectively $\mathbf{l} \vdash k_2$. The *length* of the partitions, that is, $I$, will be denoted as $|\mathbf{m}| = |\mathbf{l}|$. Let $\mathbf{l}^+$ be the cumulative sum of $\mathbf{l}$, that is, the sequence $(l_1 + \cdots + l_i)_{1 \leq i \leq I-1}$. Using $a^{\mathbf{m}}$ as a shortcut for the product of the $a^{m_i}$, the main step in the proof of (ii) above is the following inequality, valid for any positive exponents $\alpha$, $\beta$ and $\gamma$ and for $a$ small enough that $4ca^{\gamma/2} < 1$: Uniformly in $k_1$ and $k_2$,

$$(3.4) \qquad \sum_{\mathbf{m} \vdash k_1, \mathbf{l} \vdash k_2, |\mathbf{l}| = |\mathbf{m}|} a^{\alpha \mathbf{m} + \beta \mathbf{l} + \gamma \mathbf{l}^+} c^{|\mathbf{l}|} \leq Ca^{\alpha k_1/2 + \beta k_2}.$$

The direct use of this inequality and similar notation will greatly simplify the writing of the proof in the next section.

REMARK 2.   One could describe the behavior of the system in the proof of point (ii) in a different way. Let $m_t$ (resp. $l_t$) be the value of $k$ corresponding to the last $\lambda_k$ [resp. $\mathcal{C}(z, r_\tau a^k)$] discovered by $\gamma$ by time $t$ (or 0 if $t < \tau_1$, resp. $t < T_1$). The process $(m_t, l_t)$ takes values in $\mathbb{N}^2$; looking at it at times $T_{l_1 + \cdots + l_i}$, one can couple it with a discrete-time Markov chain $(M_i, L_i)$ in $\Omega = \mathbb{N}^2 \cup \{\Delta\}$, with an absorbing state $\Delta$ and transition probabilities given by the following:

- $P(M_{i+1} = M_i + m, L_{i+1} = L_i + l | M_i, L_i) = 0$ if $m \leq 0$ or $l \leq 0$,
- $P(M_{i+1} = M_i + m, L_{i+1} = L_i + l | M_i, L_i) \leq Ca^{\eta m + (8/\kappa - 1)M_i/2 + (1 - \kappa/8)l}$ if $m > 0$ and $l > 0$.

The probability estimate provided by point (ii) of the previous lemma is then bounded above by the probability that this Markov chain, started at $(0,0)$, reaches the domain $D_{k_1,k_2} = [\![k_1, \infty[\![ \times [\![k_2, \infty[\![$. Such a probability can be estimated by summing the probabilities of all possible paths going from $(0,0)$ to $U_{k_1,k_2}$ (which corresponds to the proof we just gave), or by finding an appropriate super-harmonic function on $\mathbb{N}^2$. However, we could not find a simple expression for such a super-harmonic function.

3.2. *The proof.* Applying Lemma 6, Corollary 5 and the strong Markov property, we obtain the following estimate (which we will refer to as the *main estimate*): For every good time $t$ and every radius $r \in (0, r_t)$,

$$(3.5) \qquad P(\gamma([t, \infty)) \cap \mathcal{B}(z, r) \neq \varnothing | \mathcal{F}_t, t < \infty) \leq C \left( \frac{r}{r_t} \right)^s \left( \frac{r_t}{\rho_t} \right)^{s_b/2},$$

where we define the hull and boundary exponents by

$$s = 1 - \frac{\kappa}{8} \quad \text{and} \quad s_b = \frac{8}{\kappa} - 1$$

and where $C$ depends only on $\kappa$. The way to obtain the required second-moment estimate from this upper bound is actually quite similar in spirit to the way we proved point (ii) of Lemma 8: We will split the event that $\gamma$ hits two small disks according to the order in which it visits a finite family of circles, estimate each of these individual probabilities as a product using the strong Markov property, and then sum over all possibilities. The notation is quite heavier than previously, though.

Introduce a small constant $a \in (0, 1)$ (the value of which will be determined later) and let $\delta_n = a^n \delta$. We will split the event

$$E_\varepsilon(z, z') := \{ \mathcal{B}(z, \varepsilon) \cap \gamma([0, \infty)) \neq \varnothing, \mathcal{B}(z', \varepsilon) \cap \gamma([0, \infty)) \neq \varnothing \}$$

according to the order in which the processes $(r_t)$, $(r_t')$, $(\rho_t)$ and $(\rho_t')$ reach the values $\delta_n$. For convenience, let $\bar{n} = \lfloor \log(\varepsilon/\delta)/\log a \rfloor$: It is sufficient to estimate the probability that both $(r_t)$ and $(r_t')$ reach the value $\delta_{\bar{n}}$.

Let $T_n$ (resp. $T_n'$, $\tau_n$, $\tau_n'$) be the first time when $r_t$ (resp. $r_t'$, $\rho_t$, $\rho_t'$) is not greater than $\delta_n$ (or infinity, if such a time does not exist). We will call all these stopping times *discovery times*.

LEMMA 9. *For all $n, n' > 0$, we have the following (as well as their counterparts obtained by exchanging the roles of $z$ and $z'$) if all the involved stopping times are finite:*

1. $\gamma_{\tau_n} \in \mathcal{C}(z, \delta_n) \cap \mathcal{B}_{\tau_n}(\delta_n)$; *in particular, $|\gamma_{\tau_n} - z| = \delta_n$;*
2. $T_n < T_{n+1}$ *and* $\tau_n < \tau_{n+1}$;
3. $\tilde{T}_n < \tau_n$;
4. *If $T_n < T_{n'}'$, then $\tilde{T}_n < T_{n'}'$, and similarly replacing $T$ (resp. $T'$, resp. both) by $\tau$ (resp. $\tau'$, resp. both).*

PROOF. Point (i) was proved as part of Lemma 8; point (ii) is then obvious and point (iii) is a direct consequence of Lemma 7, so only (iv) requires a proof.

Assume that $T_n < T_{n'}'$. Let $\zeta$ be a curve going from $z$ to $z'$, obtained by concatenating $\zeta_1 \subset \mathcal{B}_{T_n}(\delta_n)$, $\gamma([T_n, T_{n'}'])$ and $\zeta_2 \subset \mathcal{B}_{T_{n'}'}(\delta_{n'})$. Such a curve has to cross $\lambda'_{T_n}$ (by definition), and that can only happen on $\gamma((T_n, T'_{n'}))$ because the distance between $\mathcal{E}_{T_n}$ and $z$ (resp. $z'$) is greater than $\delta_n$ (resp. $\delta_{n'}$). This is equivalent to saying that $\tilde{T}_n < T'_{n'}$, which is what we wanted. The same reasoning applies when replacing $T$ by $\tau$ and/or $T'$ by $\tau'$. $\quad\square$

Here is a somewhat informal description of the construction we will do. Assume that $\gamma$ hits both $\mathcal{B}(z, \delta_{\bar{n}})$ and $\mathcal{B}(z', \delta_{\bar{n}})$. In order to do it, it has to cross all the circles of radii $\delta_n$, $n \leq \bar{n}$ around $z$ and $z'$, and it will do so in a certain order, coming back to the separator set $\mathcal{E}$ between explorations around $z$ and around $z'$ [this is the meaning of point (iv) of the previous lemma]. We call a *task* the time interval spanning between two successive such returns on $\mathcal{E}$. The conditional probability that a given task is performed, conditionally to its past, is then given by the main estimate (3.5), and the rest of the construction is then very similar to the proof of point (ii) in Lemma 8.

Let $S_0 = 0$ and define the stopping times $S_i$ and $S'_i$ for $i > 0$, inductively, as follows:

- $S'_i = \mathrm{Min}(\{T_n, T'_n, \tau_n, \tau'_n\} \cap (S_{i-1}, \infty))$, that is, $S'_i$ is the first discovery time after $S_{i-1}$, if such a time exists;
- $S_i = \tilde{S}'_i$ if $S'_i$ is a $T_n$ or a $\tau_n$, $S_i = \check{S}'_i$ if $S'_i$ is a $T'_n$ or a $\tau'_n$—again if such a time exists.

Continue the construction until the first ball of radius $\delta_{\bar{n}}$ is hit by the curve; let $I$ be chosen in such a way that this happens at time $S'_{I-1}$. The curve still has to come back to $\mathcal{E}$ after that, so $S_{I-1}$ is well defined. Then, simply let $S_I$ be the hitting time of the second ball of radius $\delta_{\bar{n}}$. We call *task* a time-interval of the form $(S_{i-1}, S_i]$; $I$ is then simply the number of tasks. Following our construction, the last task is different from the others and will have to be treated specially.

For each $i \leq I$, let $J_i$ (resp. $K_i$, $J'_i$, $K'_i$) be the largest integer $n > 0$ for which $\tau_n$ (resp. $T_n$, $\tau'_n$, $T'_n$) is smaller than $S_i$, if such an integer exists, and $0$ if it does not. By construction,

$$\delta_{K_i+1} < r_{S_i} \leq \delta_{K_i}, \qquad \delta_{J_i+1} < \rho_{S_i},$$

and similar inequalities hold for $r'_{S_i}$ and $\rho'_{S_i}$.

So, we obtain a sequence of quadruples $(J_i, K_i, J'_i, K'_i)_{i \leq I}$, which is not Markovian but on which we can say enough to obtain the needed second-moment estimate. First of all, for each $i < I$, at least one of $J_{i+1}$, $K_{i+1}$, $J'_{i+1}$, $K'_{i+1}$ is larger than its counterpart at index $i$; but if $J_{i+1} > J_i$ or $K_{i+1} > K_i$, then $J'_{i+1} = J'_i$ and $K'_{i+1} = K'_i$, by point (iv) of Lemma 9. The main estimate implies the following bound: for each $k > 0$,

$$(3.6) \quad \begin{aligned} &P((J_{i+1}, K_{i+1}, J'_{i+1}, K'_{i+1}) = (J_i, K_i + k, J'_i, K'_i)|\mathcal{F}_{S_i}) \\ &\leq C a^{s_b(K_i - J_i)/2} a^{sk}. \end{aligned}$$

Point (ii) of Lemma 8 then says that, for every $j > 0$ and $k \geq 0$, and if $i < I - 1$,

$$(3.7) \quad \begin{aligned} P((J_{i+1}, K_{i+1}, J'_{i+1}, K'_{i+1}) &= (J_i + j, K_i + k, J'_i, K'_i)|\mathcal{F}_{S_i}) \\ &\leq Ca^{\eta(J_i + j)}a^{sk}. \end{aligned}$$

The first of these bounds also applies in the second case, still as a consequence of the main estimate, so we get a last estimate for the last step: for every $j > 0$,

$$(3.8) \quad \begin{aligned} P((J_I, K_I, J'_I, K'_I) &= (J_{I-1} + j, \bar{n}, J'_{I-1}, \bar{n})|\mathcal{F}_{S_{I-1}}) \\ &\leq Ca^{s_b(K_{I-1} - J_{I-1})/2}a^{sk}. \end{aligned}$$

Notice that here and from now on, as the estimates on radii we get from the values of the $J$ and $K$ are only valid up to a multiplicative factor of order $a$, the constants $C$ appearing in the estimates now depend on the value of $a$.

Lemma 10 (Reduction). *With the previous notation, any jump of the second kind, that is, corresponding to equation* (3.7), *and such that both $j$ and $k$ are positive, satisfies $J_i + j = K_i$. In other words, in such a case $\gamma$ closes a fjord of width comparable to $r_{S_i}$ and then approaches $z$ before going back to the separator set (see Figure 5).*

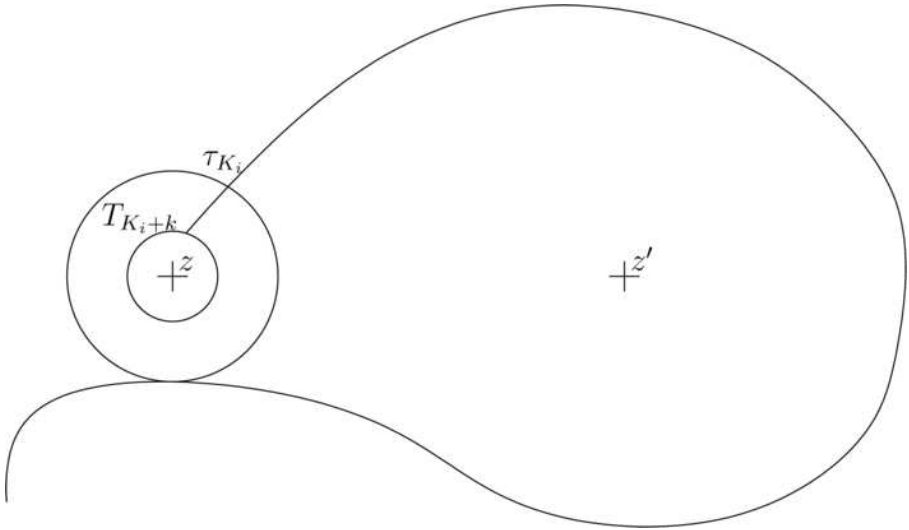

Fig. 5. *Second moments: Reduction.*

Proof. First notice that Lemma 7 ensures that $J_i + j \leq K_i$; assume that $J_i + j < K_i$. Let $\tau = \tau_{J_i+j}$ and $T = T_{K_i+k}$, and assume that both lie between $S_i$ and $S_{i+1}$ and that $\tau < T$: At time $\tau$, the situation is similar to the one in Lemma 8. Construct $\zeta'$ as in the proof of Lemma 8, and continue $\zeta'$ to a crosscut $\hat{\zeta}'$ separating $z'$ from infinity, and still contained in $\mathcal{B}_\tau(\delta_{J_i+j})$—this is possible by the definition of $\tau$. Let $\tau'$ be the first time after $\tau$ such that $\gamma_{\tau'} \in \partial \mathcal{B}_\tau(\delta_{K_i})$ (necessarily $\tau' < T$). Let $\zeta''$ be a continuous curve joining $\gamma_{\tau'}$ to $z$ inside $\mathcal{B}_\tau(\delta_{K_i})$. Because the unbounded component of $H_\tau \setminus \lambda'_\tau$ is simply connected, the concatenation of $\zeta''$ and $\gamma([\tau', \tau])$ is homotopic to $\zeta'$, and in particular, it separates $\lambda'_\tau$ from infinity. Hence, $\mathcal{B}_\tau(\delta_{K_i})$ separates $z'$ from infinity in $H_{\tau'}$, and in particular, $\tau_{K_i} < T$. This contradicts the hypothesis $J_i + j < K_i$.

The same construction applies if, among the discovery times lying between $S_i$ and $S_{i+1}$, there is a $\tau$-time before a $T$-time. The last case to consider is then when $T = T_{K_i+k}$ happens before any closing time. In particular, if such is the case, $\mathcal{B}_T(\delta_{J_i+1})$ does not separate $z'$ from infinity in $H_T$. The same proof as that of Lemma 7 then shows that $\tau_{J_i+1} > \tilde{T} = S_{i+1}$, which again is a contradiction.  □

This reduction means that, as far as reaching probabilities are concerned, the jump from $(J_i, K_i)$ to $(J_i+j, K_i+k)$ behaves exactly like the succession of two jumps, from $(J_i, K_i)$ to $(J_i+j, K_i) = (K_i, K_i)$ (which happens with probability not greater than $Ca^{\eta(J_i+j)}$) and then from $(K_i, K_i)$ to $(K_i, K_i+k)$ (which happens with probability not greater than $Ca^{sk}$): Up to replacing $C$ with $C^2$, we can assume in the computations that jumps of the type corresponding to equation (3.7) only happen with $k = 0$, and

$$(3.9) \quad P((J_{i+1}, K_{i+1}, J'_{i+1}, K'_{i+1}) = (J_i+j, K_i, J'_i, K'_i) | \mathcal{F}_{S_i}) \leq Ca^{\eta(J_i+j)}.$$

Again, this estimate does not hold for the very last task—which, according to this decomposition, would correspond to the fact that the last two jumps are a jump of $J$ followed by a jump of $K$, this last one reaching the value $\bar{n}$. In that case, with the notation in equation (3.8), we would have $j = K_{I-1} - J_{I-1}$ from the previous Lemma. Assuming that $\eta < s_b/2$, which we can do, we then obtain the following estimate for the second-to-last jump in the previous reduction (the last jump always involves $K$):

$$(3.10) \quad \begin{aligned} P((J_{I-1}, K_{I-1}, J'_{I-1}, K'_{I-1}) &= (J_{I-2}+j, K_{I-2}, J'_{I-2}, \bar{n}) | \mathcal{F}_{S_{I-2}}) \\ &\leq Ca^{\eta j}. \end{aligned}$$

In other words, here only the last jump of $J$ appears in the exponent as opposed to the end-value in the other cases. (In fact, using the main estimate here amounts to discarding that last jump of $J$ entirely, but writing the estimate this way makes for a slightly more pleasant computation below.) This turns out to be enough for our purposes.

All we need to do then is to estimate the probability that $(K_i, K_i')$ reaches $(\bar{n}, \bar{n})$. With this formulation, it would be nice to give a super-harmonic function associated to the process, but despite our best effort, we could not find such a function. The natural candidate would be of the form $C'a^{s(\bar{n}-K)+\eta(K-J)/2}$, but this might fail to be super-harmonic along the diagonal—the reason being that $C$ now depends on $a$.

So, we will apply the same strategy as in the proof of point (ii) of Lemma 8, namely, sum over all possible paths starting at $(0,0,0,0)$ and ending on $(J, \bar{n}, J', \bar{n})$ for some $J, J' \in [\![0, \bar{n}]\!]$ (which we will call *good paths*). This leads to rather unpleasant computations, but the general strategy should be clear enough.

Look first at the components $J_i$ and $K_i$ of the walk: Along a good path, the jumps of $(J_i, K_i)$ affect either its first or its second coordinate. Let $n \geq 0$ be the number of jumps affecting $(J_i)$, and let $j_1, \ldots, j_n > 0$ be their lengths. Then, for $0 \leq i \leq n$, let $l_i \geq 0$ be the number of jumps affecting $K$ between the $i$th and $(i+1)$st jumps of $J$, and let $k_{i,1}, \ldots, k_{i,l_i} > 0$ be their lengths— with the obvious abuse of notation that $l_0$ is the number of jumps of $K$ before the first jump of $J$, and $l_n$ is the number of jumps of $K$ after the last jump of $J$. In particular, the sum of all $k_{i,j}$ is equal to $\bar{n}$. Define the integers $n'$, $j_i'$, $l_i'$ and $k_{i,j}'$ in a similar fashion to describe the behavior of $(J_i', K_i')$.

Notice that, just before the jump corresponding to $j_i$, the value of $J$ is $j_1 + \cdots + j_{i-1}$. Besides, let $d_i \geq 0$ be the value of $K - J$ just after that jump ($d_i$ is the difference between the sum of the $k$'s and that of the $j$'s so far): Just before the jump corresponding to $k_{i,j}$, the value of $K - J$ is equal to $d_i + k_{i,1} + \cdots + k_{i,j-1}$. This is sufficient to estimate the probability that a given path occurs: It will be given by the product along the path of the (conditional, given the past) probabilities of the individual steps, which is bounded above by the product of two terms, namely,

$$
\begin{aligned}
A := {}& a^{-\eta(j_1 + \cdots + j_{n-1})} \\
& \times \prod_{j=1}^{l_0} C(a^{k_{0,j}})^s (a^{k_{0,1} + \cdots + k_{0,j-1}})^{s_b/2} \\
& \times \prod_{i=1}^{n} \left[ C(a^{j_1 + \cdots + j_i})^\eta \prod_{j=1}^{l_i} C(a^{k_{i,j}})^s (a^{d_i + k_{i,1} + \cdots + k_{i,j-1}})^{s_b/2} \right].
\end{aligned}
$$

Here, the empty products are equal to 1 by convention. The term $a^{-\eta(j_1 + \cdots + j_{n-1})}$ accounts for the difference in the very last task, where as was pointed out above, only the last jump of $J$, that is, $j_n$, appears in the exponent. That task might be a jump toward $z'$, or not involve a jump of $J$ at all, in which case the factor would not be needed, but having it in all cases makes the computation more symmetric—it not smaller than 1 anyway.

Also define the corresponding $A'$ involving $n'$, $j'_i$, $l'_i$ and $k'_{i,j}$. With the shortcut notation used in the previous subsection, letting $\mathbf{j} = (j_i)$ and $\mathbf{k}_i = (k_{i,j})$, so that $|\mathbf{j}| = n$ and $|\mathbf{k}_i| = l_i$, this becomes

$$A = a^{-\eta \mathbf{j}_{n-1}^+} C^{|\mathbf{j}| + \sum |\mathbf{k}_i|} a^{s \mathbf{k_0}} a^{s_b \mathbf{k}_0^+/2} a^{\eta \mathbf{j}^+} \prod_{i=1}^{n} [a^{s \mathbf{k}_i} a^{s_b \mathbf{k}_i^+/2} a^{|\mathbf{k}_i| d_i s_b / 2}].$$

Rewriting the product using the fact that $\sum k_{i,j} = \bar{n}$, and letting $j_0 = d_0 = 0$ for ease of notation, we obtain

$$A = a^{-\eta \mathbf{j}_{n-1}^+} C^{n + \sum l_i} a^{s \bar{n}} \prod_{j=1}^{l_0} (a^{k_{0,1} + \cdots + k_{0,j-1}})^{s_b/2}$$

$$\times \prod_{i=1}^{n} \left[ (a^{j_1 + \cdots + j_i})^{\eta} \prod_{j=1}^{l_i} (a^{d_i + k_{i,1} + \cdots + k_{i,j-1}})^{s_b/2} \right]$$

$$= a^{-\eta \mathbf{j}_{n-1}^+} C^{n + \sum l_i} a^{s \bar{n}} \prod_{i=0}^{n} \left[ (a^{j_1 + \cdots + j_i})^{\eta} \prod_{j=1}^{l_i} (a^{d_i + k_{i,1} + \cdots + k_{i,j-1}})^{s_b/2} \right]$$

$$= C^{n + \sum l_i} a^{s \bar{n}} a^{\eta j_n} \prod_{i=0}^{n} \left[ a^{(n-i)\eta j_i + l_i d_i s_b / 2} \prod_{j=1}^{l_i - 1} a^{(l_i - j) k_{i,j} s_b / 2} \right].$$

Indeed, each term $a^{j_i}$ for $i < n$ appears $n - i$ times in the product, the term $a^{j_n}$ appears once [in other terms, $a^{j_i}$ appears $(n - i) \vee 1$ times], the term $a^{d_i}$ appears $l_i$ times and the term $a^{k_{i,j}}$ appears $l_i - j$ times.

We still have to sum the product $AA'$ over all the good paths. Notice first that giving the values of $n$, $j_i$, $l_i$ and $k_{i,j}$, $n'$, $j'_i$, $l'_i$ and $k'_{i,j}$ is not sufficient to specify the path of $(J, K, J', K')$, because it says nothing about the way the jumps of $(J, K)$ and $(J', K')$ are intertwined; however, there are at most

$$\binom{n + \sum l_i + n' + \sum l'_i}{n + \sum l_i} \leq 2^{n + \sum l_i} 2^{n' + \sum l'_i}$$

such intertwinings. Hence, up to doubling of the constant $C$, it is sufficient to sum $AA'$ over the values of $n$, $j_i$, $l_i$ and $k_{i,j}$, $n'$, $j'_i$, $l'_i$ and $k'_{i,j}$. The sum will factor into two terms, one involving the terms around $z$ and the other the terms around $z'$, and these two factors are equal. Hence, an upper bound of the probability that $(K, K')$ reaches $(\bar{n}, \bar{n})$ is given by $B^2$, where

$$B := \sum_{\mathbf{j}, \mathbf{k}_i} a^{-\eta \mathbf{j}_{n-1}^+} C^{|\mathbf{j}| + \sum |\mathbf{k}_i|} a^{s \mathbf{k_0}} a^{s_b \mathbf{k}_0^+/2} a^{\eta \mathbf{j}^+} \prod_{i=1}^{n} [a^{s \mathbf{k}_i} a^{s_b \mathbf{k}_i^+/2} a^{|\mathbf{k}_i| d_i s_b / 2}],$$

with a sum taken over all values of the indices leading to a good path.

First, let $k_i = \sum_{j=1}^{l_j} k_{i,j}$ (and notice that $d_i = d_{i-1} + k_{i-1} - j_i$). We can rewrite $B^2$ as

$$B^2 = a^{2s\bar{n}} \left[ \sum_{n=0}^{\infty} C^n \prod_{i=0}^{n} \left[ \sum_{j_i, k_i} a^{[(n-i)\vee 1]\eta j_i} \sum_{l_i, k_{i,j}} C^{l_i} a^{l_i d_i s_b/2} \prod_{j=1}^{l_i-1} a^{(l_i-j)k_{i,j}s_b/2} \right] \right]^2,$$

where the innermost sum is taken over all choices of $l_i$ and $k_{i,j}$ satisfying $\sum k_{i,j} = k_i$, and where the sum over $j_i$ is in fact not present for $i = 0$. This in turn can be considered as a sum over the $k_{i,j}$ for $j \leq l_i - 1$ with sum smaller than $k_i$. The case $k_i = l_i = 0$ needs to be treated separately here, and we get

$$B^2 \leq a^{2s\bar{n}} \left[ \sum_{n=0}^{\infty} C^n \prod_{i=0}^{n} \left[ \sum_{j_i} a^{[(n-i)\vee 1]\eta j_i} \left( 1 + \sum_{k_i, l_i > 0} C^{l_i} a^{l_i d_i s_b/2} \right. \right. \right.$$
$$\left. \left. \left. \times \prod_{j=1}^{l_i-1} \left( \sum_{k>0} a^{(l_i-j)ks_b/2} \right) \right) \right] \right]^2.$$

The sum over $k$ can be computed explicitly, it is convergent because $j < l_i$ and its value is smaller than $2a^{(l_i-j)s_b/2}$ if $a$ is chosen small enough, which we will assume from now on. The product over $j$ is then equal to $2^{l_i-1}a^{l_i(l_i-1)s_b/4}$. The terms with an exponent linear in $l_i$ can be factored into $C^{l_i}$—note that $C$ now depends on $a$—leading to

$$B^2 \leq a^{2s\bar{n}} \left[ \sum_{n=0}^{\infty} C^n \prod_{i=0}^{n} \left[ \sum_{j_i} a^{[(n-i)\vee 1]\eta j_i} \left( 1 + \sum_{d_i, l_i} C^{l_i} a^{l_i d_i s_b/2} a^{l_i^2 s_b/4} \right) \right] \right]^2$$
$$\leq a^{2s\bar{n}} \left[ \sum_{n=0}^{\infty} C^n \prod_{i=0}^{n} \left[ \sum_{j_i} a^{[(n-i)\vee 1]\eta j_i} \left( 1 + \sum_{d_i \geq 0} a^{d_i s_b/2} \sum_{l_i > 0} C^{l_i} a^{l_i^2 s_b/4} \right) \right] \right]^2.$$

The sums over $l_i$ and $d_i$ are convergent, because $a < 1$, so the whole term in parentheses is bounded by a constant depending only on $\kappa$ and $a$; since it appears $n + 1$ times, up to another change in the value of $C$, we get

$$B^2 \leq a^{2s\bar{n}} \left[ \sum_{n=0}^{\infty} C^{n+1} \prod_{i=1}^{n} \left[ \sum_{j_i > 0} a^{[(n-i)\vee 1]\eta j_i} \right] \right]^2.$$

Summing over all values of $j_i > 0$, we obtain (if $a$ is small enough)

$$B^2 \leq a^{2s\bar{n}} \left[ \sum_{n=0}^{\infty} C^{n+1} 2a^{\eta} \prod_{i=1}^{n-1} \left( 2a^{(n-i)\eta} \right) \right]^2.$$

Up to yet another increase of $C$, the factor 2 can be made part of it, and the product over $i$ can be computed explicitly:

$$B^2 \leq a^{2s\bar{n}} \left[ \sum_{n=0}^{\infty} C^{n+1} a^{\eta} a^{n(n-1)\eta/2} \right]^2.$$

This last sum is again convergent, so we obtain $B^2 \leq Ca^{2s\bar{n}}$, with a constant $C$ depending only on $\kappa$ and $a$.

Putting everything together, assuming $E_\varepsilon(z, z')$ holds, first $\gamma$ has to reach $\mathcal{E}$, and this happens with probability of order $\delta^s$. Then, conditionally to the process up to this hitting time, we can apply the previous reasoning which says that the conditional probability to hit both disks of radius $\delta_{\bar{n}}$ is bounded above by $Ca^{2\bar{n}s}$ where $C$ depends only on $a$ and $\kappa$. Hence,

$$P(E_\varepsilon(z, z')) \leq C\delta^s a^{2\bar{n}s}.$$

Notice that $a^{\bar{n}} \leq (\varepsilon/\delta)a^{-1}$ to finally obtain

$$P(E_\varepsilon(z, z')) \leq Ca^{-2s}\frac{\varepsilon^{2s}}{\delta^s},$$

which is precisely the estimate we were looking for.

**4. The occupation density measure.** As a side remark, let us consider the proof of the lower bound for the dimension (cf. Section 1). It is based on the construction of a Frostman measure $\mu$ supported on the path, constructed as a subsequential limit of the family $(\mu_\varepsilon)$ defined by their densities with respect to the Lebesgue measure on the upper-half plane:

$$d\mu_\varepsilon(z) = \varepsilon^{-s}\mathbb{1}_{z \in C_\varepsilon}|dz|.$$

Then, $\mu$ is a random measure with correlations between $\mu(A)$ and $\mu(B)$, for disjoint compact sets $A$ and $B$, decaying as a power of their inverse distance. So, at least formally, it behaves in this respect like a conformal field: the one-point function (corresponding to the density of $\mu$) is not welldefined, because $\mu$ is singular to the Lebesgue measure, but the two-point correlation

$$\lim_{\delta \to 0} \delta^{-4} \operatorname{Cov}(\mu(\mathcal{B}(z, \delta)), \mu(\mathcal{B}(z', \delta)))$$

behaves like $d(z, z')^{-1+\kappa/8}$.

A little more can be said about this measure, or about its expectation. The proof of the estimate for $P(\gamma \cap \mathcal{B}(z, \varepsilon) \neq \varnothing)$ can be refined in the following way: When we apply the stopping theorem (2.8), saying that the diffusion conditioned to survive has a limiting distribution shows that

$$E\left[\sin\left(\frac{\alpha_s}{2}\right)^{8/\kappa-1}\Big| S \geq s\right]$$

has a limit $\lambda$ when $s \to \infty$, and that this limit depends only on $\kappa$. So what we get out of the construction in Section 2 is

$$P\left(\exists t > 0 : |g'_t(z)| \geq \frac{\Im(z)}{\varepsilon}\right) \underset{\varepsilon \to 0}{\sim} \lambda(\kappa)\left(\frac{\varepsilon}{\Im z}\right)^{1-\kappa/8}(\sin(\arg(z)))^{8/\kappa-1}.$$

This lead us to an estimate on $P(d(z,\gamma) < \varepsilon)$ by the Köbe 1/4 theorem; but it is also natural to measure the distance to $\gamma$ by the modulus of $g'$. We can now define

$$\phi_1(z) = \lim_{\varepsilon \to 0} \varepsilon^{\kappa/8 - 1} P\left(\exists t > 0 \colon |g'_t(z)| \geq \frac{\Im(z)}{\varepsilon}\right) \colon$$

the previous estimate boils down to

$$\phi_1(z) = \lambda(\kappa) \Im(z)^{\kappa/8 - 1} \sin(\arg z)^{8/\kappa - 1},$$

and by the construction of $\mu$, we obtain that, for every Borel subset $A$ of the upper-half plane,

$$E(\mu(A)) \asymp \int_A \phi_1(z) |\, dz|,$$

with universal constants.

It is then possible to do this construction for several points; note first that the second moment estimate can actually be written as

$$P(\{z, z'\} \subset C_\varepsilon) \asymp \frac{\varepsilon^{2(1 - \kappa/8)}}{|z - z'|^{1 - \kappa/8} \Im((z + z')/2)^{1 - \kappa/8}},$$

as long as both $\Im(z)$ and $\Im(z')$ are bounded below by $|z - z'|/M$ for some fixed $M > 0$. Indeed, the upper bound is exactly what we derived in the previous section, and the lower bound is provided by the term $n = n' = 0$ in the sum. Hence, any subsequential limit $\psi(z, z')$, as $\varepsilon$ vanishes, of

$$\varepsilon^{2(\kappa/8 - 1)} P(\{z, z'\} \subset C_\varepsilon)$$

satisfies $\psi(z, z') \asymp \phi_2(z, z')$ for some fixed function $\phi_2$, with constants depending only on $\kappa$. The second moment estimate then shows that

$$\phi_2(z, z') \underset{z' \to z}{\asymp} \frac{\phi_1(z)}{|z - z'|^{1 - \kappa/8}},$$

that is, $\phi_2$ behaves like a correlation function when $z$ and $z'$ are close to each other.

The general case of $n$ points, $n \geq 2$, can be treated in the same fashion. First, the derivation of second moments admits a generalization to $n$ points, as follows. Let $(z_i)_{1 \leq i \leq n}$ be $n$ distinct points in $\mathbb{H}$, such that their imaginary parts are large enough (bigger than, say, $18n$ times the maximal distance between any two of them). We use them to construct a Voronoi tessellation of the plane; denote by $C_i$ the face containing $z_i$, and by $\delta_i$ the (Euclidean) distance between $z_i$ and $\partial C_i$. Let $\mathcal{C}(z_0, \delta_0)$ be the smallest circle containing all the discs $\mathcal{B}(z_i, \delta_i)$. Last, let $\mathcal{E}$ be the "separator set" between the $z_i$'s, defined as

$$\mathcal{E} = \mathcal{C}(z_0, \delta_0) \cup \left[\left(\bigcup_{i=1}^{n} \partial C_i\right) \cap \mathcal{B}(z_0, \delta_0)\right].$$

It is the same as defined previously in the case $n = 2$.

The previous proof can then be adapted to show that

$$P(\{z_1, \ldots, z_n\} \subset C_\varepsilon) \asymp \left( \frac{\delta_0 \varepsilon^n}{\prod \delta_i} \right)^{1 - \kappa/8}$$

(using radii $\delta_i a^k$ for the circles around $z_i$). In the case $n = 2$, we have $\delta_1 = \delta_2 = \delta_0/2$, so this estimate is exactly the same as previously. So, it makes sense to take a (subsequential) limit, as $\varepsilon$ tends to 0, of

$$\varepsilon^{n(\kappa/8 - 1)} P(\{z_1, \ldots, z_n\} \subset C_\varepsilon),$$

and all possible subsequential limits are comparable to a fixed symmetric function $\phi_n$.

The behavior of $\phi_n(z_1, \ldots, z_n)$ when $z_n$ approaches the boundary is then given by the boundary term in Proposition 4, that is, $\phi$ behaves like $(\Im z_n)^{8/\kappa - 1}$ there. Last, it is easy to see that, when $z_n$ tends to $z_1$, $\phi_n(z_1, \ldots, z_n)$ has a singularity which is comparable to $|z_n - z_1|^{\kappa/8 - 1}$; in other words, we have a recursive relation between all the $\phi_n$'s, given by

$$(4.1) \qquad \phi_n(z_1, \ldots, z_n) \underset{z_n \to z_1}{\asymp} \frac{\phi_{n-1}(z_1, \ldots, z_{n-1})}{|z_n - z_1|^{1 - \kappa/8}},$$

$$(4.2) \qquad \phi_n(z_1, \ldots, z_n) \underset{\Im z_n \to 0}{\asymp} \phi_{n-1}(z_1, \ldots, z_{n-1}) \cdot (\Im z_n)^{8/\kappa - 1}.$$

These relations are very similar to some of those satisfied by the correlation functions in conformal field theory. In fact, it is possible to push the relation further, in two ways. First, we can look at the evolution of the system in time. This corresponds to mapping the whole picture by the map $g_t - \beta_t$, and this map acts on the discs of small radius around the $z_i$'s like a multiplication of factor $|g_t'(z_i)|$ (as long as $K_t$ remains far away from the $z_i$'s, which we may assume if $t$ is small enough). Hence, the process

$$Y_t^n := \left( \prod |g_t'(z_i)|^{1 - \kappa/8} \right) \phi_n(g_t(z_1) - \beta_t, \ldots, g_t(z_n) - \beta_t)$$

(defined as long as all the $z_i$'s remain outside $K_t$) is a local martingale. We can apply Itô's formula to compute $dY_t^n$, and write that the drift term has to be 0 at time 0 to obtain a PDE satisfied by $\phi_n$.

Note though that the formula involves the modulus of $g_t'$, meaning that the equation we would obtain cannot be expressed in terms of complex derivatives of $g_t$ only, and that we have to introduce derivatives with respect to the coordinates. This is also the case for the second-order term in Itô's formula: Since $\beta$ is a real process, we would obtain terms involving second derivatives of $\phi_n$ with respect to the $x$-coordinates of the arguments. To sum it up, it would be an ugly formula without the correct formalism—which is why we do not put it here. The formula is much nicer when considering points on the boundary of the domain—compare [5].

The last thing we can do is study what happens if we add one point $z_{n+1}$ to the picture. This will add one multiplicative factor, corresponding (at least intuitively) to the conditional probability to hit $z_{n+1}$ knowing that we touch the first $n$ points already. In the case $\kappa = 8/3$ and for points on the boundary of the domain, this can be computed using the restriction property, and it leads to Ward's equations (cf. [5]). In the "bulk" (i.e., for points inside the domain), or for other values of $\kappa$, it is not clear yet how to do it.

**5. The boundary.** A natural question is the determination of the dimension of the boundary of $K_t$ for some fixed $t$, in the case $\kappa > 4$. The conjectured value is

$$\dim_H(\partial K_t) = 1 + \frac{2}{\kappa},$$

and this can now be proved for a few values of $\kappa$ for which the boundary of $K$ can be related to the path of an $\text{SLE}_{\kappa'}$ with $\kappa' = 16/\kappa$. In fact, this relation is only known in the cases where convergence of a discrete model to SLE is known, namely:

- $\kappa = 6$, where actually both $\partial K_t$ and the path of the $\text{SLE}_{\kappa'}$ are closely related to the Brownian frontier. Hence, we obtain a third derivation of the dimension of the Brownian frontier, this time through $\text{SLE}_{8/3}$.
- $\kappa = 8$: Here, $\text{SLE}_8$ is known to be the scaling limit of the uniform Peano curve and $\text{SLE}_2$ that of the loop-erased random walk (cf. [12]). Since these two discrete objects are closely related through Wilson's algorithm, this shows that the local structure of the $\text{SLE}_2$ curve and the $\text{SLE}_8$ boundary are the same, and in particular, they have the same dimension.

So we obtain one additional result here:

COROLLARY 11. *Let* $(K_t)$ *be a chordal* $\text{SLE}_8$ *in the upper-half plane: Then, for all* $t > 0$, *the boundary of* $K_t$ *almost surely has Hausdorff dimension* 5/4.

It would be nice to have a direct derivation of the general result, without using the "duality" between $\text{SLE}_\kappa$ and $\text{SLE}_{16/\kappa}$; but it is not even clear how to obtain a precise estimate of the probability that a given ball intersects the boundary of $K_1$.

**Acknowledgments.** Part of this work was carried out during my stay at the Mittag–Leffler institute whose hospitality and support is acknowledged. I also thank Peter Jones, Greg Lawler and Wendelin Werner for very useful discussions, and I especially wish to thank Oded Schramm for his help and patience in reviewing several draft versions of this paper.

UMPA—ENS LYON—CNRS UMR 5669
46 ALLÉE D'ITALIE
F-69364 LYON CEDEX 07
FRANCE
E-MAIL: vbeffara@ens-lyon.fr